\documentclass[12pt,reqno]{article}
\usepackage{amsmath,amsthm,amsfonts,amssymb,amscd}
\usepackage[dvips]{graphicx}
\usepackage{psfrag}

\input{psfig}
\theoremstyle{plain}
\newtheorem{thm}{Theorem}[section]

\newtheorem{prop}[thm]{Proposition}

\newtheorem{lemma}[thm]{Lemma}

\newtheorem{defi}[thm]{Definition}

\newtheorem{maintheorem}{Theorem}

\title{Sufficient conditions for robustness of attractors}
\author{C. A. Morales, M. J. Pacifico
\thanks{
Partially supported by CNPq, FAPERJ and PRONEX/DYN. SYS.}}

\begin{document}
\maketitle
\begin{abstract}
A recent problem in dynamics is to determinate
whether an attractor $\Lambda$ of a $C^r$ flow
$X$ is $C^r$ robust transitive or not.
By {\em attractor}
we mean a transitive set
to which all positive orbits close to it converge.
An attractor is $C^r$ robust transitive (or {\em $C^r$ robust} for short)
if it exhibits a neighborhood $U$
such that the set $\cap_{t>0}Y_t(U)$
is transitive for every flow $Y$ $C^r$ close to $X$.
We give sufficient conditions for robustness of attractors based on the following definitions.
An attractor
is {\em singular-hyperbolic} if
it has singularities (all hyperbolic) and is
partially hyperbolic with volume expanding central
direction \cite{MPP}.
An attractor is
{\em $C^r$ critically-robust} if it exhibits
a neighborhood $U$ such that
$\cap_{t>0}Y_t(U)$ is in the
closure of the closed orbits
is every flow $Y$ $C^r$ close to $X$.
We show that on compact $3$-manifolds
all $C^r$ critically-robust
singular-hyperbolic attractors
with only one singularity are $C^r$ robust.
\end{abstract}

\section{Introduction}

A recent problem in dynamics is to determinate
whether an attractor $\Lambda$ of a $C^r$ flow
$X$ is $C^r$ robust transitive or not.
By {\em attractor}
we mean a transitive set
to which all positive orbits close to it converge.
An attractor is $C^r$ robust transitive (or
{\em $C^r$ robust} for short)
if it exhibits a neighborhood $U$
such that the set $\cap_{t>0}Y_t(U)$
is transitive for every flow $Y$ $C^r$ close to $X$.
We give sufficient conditions for robustness of attractors based on the following definitions.
An attractor
is {\em singular-hyperbolic} if
it has singularities (all hyperbolic) and is
partially hyperbolic with volume expanding central
direction \cite{MPP}.
An attractor is
{\em $C^r$ critically-robust} if it exhibits
a neighborhood $U$ such that
$\cap_{t>0}Y_t(U)$ is in the
closure of the closed orbits
is every flow $Y$ $C^r$ close to $X$.
We show that on compact $3$-manifolds
all $C^r$ critically-robust
singular-hyperbolic attractors
with only one singularity are $C^r$ robust.
Let us state our result in a precise way.

Hereafter $X_t$ is a flow induced by a $C^r$ vector field $X$
on a compact $3$-manifold $M$.
The $\omega$-limit set
of $p\in M$ is the accumulation point set
$\omega_X(p)$ of the positive orbit of $p$.
An invariant set is {\em transitive} if it is
$\omega_X(p)$ for some point $p$ on it.
An {\em attracting set}
is a compact set realizing as
as $\cap_{t>0}X_t(U)$ for some
neighborhood $U$ and an {\em attractor} is a
transitive attracting set.
See \cite{Mi} where several definitions of attractor
are discussed.
The central definition of this paper
is the following.

\begin{defi}
\label{DEF0}
An attractor of a $C^r$ flow $X$ is
$C^r$ robust transitive (or {\em $C^r$ robust} for
short) if it exhibits a neighborhood
$U$ such that $\cap_{t>0}Y_t(U)$
is a transitive set
of $Y$ for every flow $Y$ $C^r$ close to $X$.
\end{defi}

Very recently \cite{B,P}
introduce the problem of finding
sufficient conditions for robustness of attractors.
To attact this problem we introduce the following definitions.
A compact invariant set
$\Lambda$ of $X$ is {\em partially hyperbolic}
if there are an invariant splitting
$T\Lambda=E^s\oplus E^c$
and positive constants $K,\lambda$ such that:

\begin{enumerate}
\item
{\em $E^s$ is contracting}, namely
$$
\mid\mid DX_t/E^s_x\mid \mid\leq
K e^{-\lambda t},
\,\,\,\,\,\,\forall x\in \Lambda,\,\,\forall
t>0.
$$
\item
{\em $E^s$ dominates $E^c$}, namely
$$
\mid \mid DX_t/E^s_x\mid\mid
\cdot
\mid\mid
DX_{-t}/E^c_{X_{t}(x)}
\mid\mid
\leq K e^{-\lambda t},
\,\,\,\,\,\,\forall x\in \Lambda,\,\,\forall
t>0.
$$
\end{enumerate}

The central direction $E^c$ of
$\Lambda$ is said to be {\em volume expanding}
if the additional condition
$$
\mid J(DX_t/E^c_x)\mid
\geq K e^{\lambda t}
$$
holds $\forall x\in \Lambda$, $\forall t>0$
where
$J(\cdot)$ means the jacobian.

\begin{defi}(\cite{MPP})
\label{DEF1}
An attractor is {\em singular-hyperbolic} if
it has singularities (all hyperbolic) and is partially hyperbolic
with volume expanding central direction.
\end{defi}

The most representative example
of a $C^r$ robust singular-hyperbolic
attractor is the geometric Lorenz attractor
\cite{GW}. The main result in \cite{MPP} claims
that $C^1$ robust non-trivial attractors
on compact $3$-manifolds are singular-hyperbolic.
The converse is false, namely
there are singular-hyperbolic attractors on compact $3$-manifolds which are not $C^r$ robust \cite{MPu}.
The following definition gives a further
sufficient condition for robustness.

\begin{defi}
\label{DEF2}
An attractor of a $C^r$ flow $X$ is
{\em $C^r$ critically-robust}
if there is a neighborhood $U$ of it such that
$\cap_{t>0}Y_t(U)$ is in the
closure of the closed orbits of $Y$
is every flow $Y$ $C^r$ close to $X$.
\end{defi}

Hyperbolic attractors on compact
manifolds are $C^r$ robust and
$C^r$ critically-robust for all $r$.
The geometric Lorenz attractor
\cite{GW} is
an example of a singular-hyperbolic attractor
with only one singularity
which is also $C^r$ robust and $C^r$ critically-robust.
In general singular-hyperbolic attractors
with only one singularity
may be neither $C^r$ robust nor $C^r$ critically-robust
\cite{MPu}.
Nevertheless we shall prove
that on compact $3$-manifolds
$C^r$ critically-robustness implies
$C^r$ robustness among
singular-hyperbolic attractors with
only one singularity.
More precisely one has the following.

\begin{maintheorem}
\label{thA}
$C^r$ critically-robust
singular-hyperbolic attractors
with only one singularity
on compact $3$-manifolds
are $C^r$ robust.
\end{maintheorem}

This theorem gives explicit sufficient conditions
for robustness of attractors {\em depending on the perturbed flow}.
E. Pujals is interested in
conditions {\em depending
on the unperturbed flow only}.
It would be also interesting
to determinate
whether the conclusion of Theorem \ref{thA} holds
interchanging the roles of robust
and critically-robust in the statement.
The proof of Theorem \ref{thA} relies on the results
in the recent work \cite{MP2}.
We reproduce these results
in Section 2 for the sake of completeness.
The proof of Theorem \ref{thA} is in Section 3.

\section{Singular-hyperbolic attracting sets}

In this section we
describe the results in \cite{MP2}.
Some proofs will be omitted.
We refer to \cite{MP2} for further details.
Hereafter $X$ is a $C^r$ flow
on a closed $3$-manifold $M$.
The closure of $B$ will be denoted by $Cl(B)$.
If $A$ is a compact invariant set
of $X$ we denote $Sing_X(A)$
the set of singularites of $X$ in $A$.
We denote by
$Per_X(A)$ the union of the
periodic orbits of $X$ in $A$.
A compact invariant set $H$ of
$X$ is {\em hyperbolic}
if the tangent bundle over $H$ has an
invariant decomposition
$E^s\oplus E^X\oplus E^u$
such that $E^s$ is contracting, $E^u$ is expanding
and $E^X$ is generated by the direction of $X$
\cite{PT}.
The Stable Manifold Theory \cite{HPS}
asserts the existence of the stable manifold
$W^s_X(p)$ and the unstable manifold $W^u_X(p)$
associated to $p\in H$.
These manifolds are respectively tangent to the
subspaces
$E^s_p\oplus E^X_p$ and $E^X_p\oplus E^u_p$
of $T_pM$.
In particular, $W^s_X(p)$ and $W^u_X(p)$
are well defined if $p$ belongs to a hyperbolic
periodic orbit of $X$.
If $O$ is an orbit of $X$ we denote
$W^s_X(O)=W^s_X(p)$ and $W^u_X(O)=W^u_X(p)$ for some
$p\in O$.
We observe that $W^{s(u)}_X(O)$ does not depend
on $p\in O$.
When $dim(E^s)=dim(E^u)=1$ we say that
$H$ is {\em saddle-type}.
In this case $W^s_X(p)$ and $W^u_X(p)$ are
two-dimensional submanifolds of $M$.
The maps
$p\in H\to W^s_X(p)$ and $p\in H\to W^u_X(p)$ are
continuous (in compact parts).
On the other hand,
a compact, singular, invariant set
$\Lambda$ of $X$ is {\em
singular-hyperbolic}
if all its singularities are hyperbolic and
the tangent bundle over $\Lambda$ has an
invariant decomposition
$E^s\oplus E^c$
such that $E^s$ is contracting, $E^s$
dominates $E^c$ and
$E^c$ is volume expanding
(i.e.
the jacobian of
$DX_t/E^c$ grows exponentially
as $t\to \infty$).
Again the Stable Manifold Theory
asserts the existence of the strong stable manifold
$W^{ss}_X(p)$
associated to $p\in \Lambda$.
This manifold is tangent to the subspace
$E^s_p$ of $T_pM$.
For all
$p\in \Lambda$ we
define
$W^s_X(p)=\cup_{t\in I\!\! R}W^{ss}_X(X_t(p))$.
If $p$ is regular (i.e. $X(p)\neq 0$)
then $W^s_X(p)$ is
a well defined two-dimensional submanifold
of $M$.
The map $p\in \Lambda\to W^s_X(p)$
is continuous (in compact parts)
at the regular points $p$ of
$\Lambda$.
A singularity $\sigma$ of $X$ is
{\em Lorenz-like} if its eigenvalues
$\lambda_1,\lambda_2,\lambda_3$
are real and satisfy
$$
\lambda_2<\lambda_3<0<-\lambda_3<\lambda_1
$$
up to some order.
A Lorenz-like singularity $\sigma$ is hyperbolic,
and so, $W^s_X(\sigma)$ and $W^u_X(\sigma)$
do exist.
Moreover, the eigenspace
of $\lambda_2$ is tangent
to a one-dimensional invariant manifold $W^{ss}_X(\sigma)$.
This manifold is called the {\em strong stable
manifold} of $\sigma$.
Clearly $W^{ss}_X(\sigma)$ splits
$W^s_X(\sigma)$ in two connected components.
We denote by
$W^{s,+}_X(\sigma), W^{s,-}_X(\sigma)$
the two connected components
of $W^s_X(\sigma) \setminus
W^{ss}_X(\sigma)$.

When $\Lambda$ is a singular-hyperbolic set
with dense periodic orbits we know that
every $\sigma\in Sing_X(\Lambda)$ is Lorenz-like
and satisfies
$\Lambda\cap W^{ss}_X(\sigma)=\{\sigma\}$
(see \cite{MPP}).
It follows also from
\cite{MPP} that
any compact invariant set without singularities
of $\Lambda$
is hyperbolic saddle-type.
When the periodic orbits are dense in $\Lambda$
attracting we have that
for every $p\in Per_X(\Lambda)$ there is
$\sigma\in Sing_X(\Lambda)$ such that
$$
W^u_X(p)\cap W^s_X(\sigma)\neq \emptyset.
$$
This fact is proved using the methods in
\cite{MP1}.

If $\Lambda$ is a singular-hyperbolic set
and $\sigma\in Sing_X(\Lambda)$ is Lorenz-like
we define

\begin{itemize}
\item
$
P^+=\{p\in Per_X(\Lambda):
W^u_X(p)\cap W^{s,+}_X(\sigma)\neq\emptyset\}$.
\item
$
P^-=\{p\in Per_X(\Lambda):
W^u_X(p)\cap W^{s,-}_X(\sigma)\neq\emptyset\}$.
\item
$
H^+_X=Cl(P^+)$.
\item
$
H^-=Cl(P^-)$.
\end{itemize}

These sets will play important role in the sequel.

\begin{lemma}
\label{cara}
Let $\Lambda$ be a connected, singular-hyperbolic,
attracting set with dense periodic orbits
and only one singularity $\sigma$.
Then, $
\Lambda=H^+\cup
H^-$.
\end{lemma}


Next we state a technical lemma to be used later.
If  $S$ is a submanifold we shall denote
by $T_xS$ the tangent space at $x\in S$.
A {\em cross-section} of $X$ is
a compact submanifold $\Sigma$
transverse to $X$ and diffeomorphic to the
two-dimensional square $[0,1]^2$.
If $\Lambda$ is a singular-hyperbolic set
of $X$ and $x\in \Sigma\cap \Lambda$,
then $x$ is regular and so
$W^s_X(x)$ is a two-dimensional submanifold
transverse to $\Sigma$.
In this case we define
by $W^s_X(x,\Sigma)$ the
connected component of $W^s_X(x)\cap \Sigma$
containing $x$.
We shall be interested in a
special cross-section described as follows.
Let $\Lambda$ be a singular-hyperbolic set
of a three-dimensional flow $X$
and $\sigma\in Sing_X(\Lambda)$.
Suppose that the closed orbits contained in $\Lambda$ are dense in $\Lambda$.
Then $\sigma$ is Lorenz-like
\cite{MPP}, and so, it is possible to describe the flow
using the Grobman-Hartman Theorem
\cite{dMP}.
Indeed, we can assume that
the flow of $X$ around
$\sigma$ is
the linear flow
$
\lambda_1\partial_{x_1}+\lambda_2\partial_{x_2}+
\lambda_3\partial_{x_3}
$
in a suitable coordinate system
$(x_1,x_2,x_3)\in [-1,1]^3$ around
$\sigma=(0,0,0)$.
A cross-section $\Sigma$ of $X$ is
{\em singular}
if $\Sigma$ corresponds to the submanifolds
$\Sigma^+=\{z=1\}$ or
$\Sigma^-=\{z=-1\}$ in the coordinate system
$(x_1,x_2,x_3)$.
We denote by $l^+$ and $l^-$
the curves in $\Sigma^+, \Sigma^-$
intersecting to $\{y=0\}$.
Note that $l^+,l^-$ are
contained in $W^{s,+}_X(\sigma),W^{s,-}_X(\sigma)$ respectively.
We state without proof
the following straighforward lemma.

\begin{lemma}
\label{sigma}
Let $\Lambda$ a singular-hyperbolic set
with dense periodic orbits of a three-dimensional flow $X$ and $\sigma\in Sing_X(\Lambda)$
be fixed. Then,
there are singular
cross-sections $\Sigma^+,\Sigma^-$
as above such that
every orbit of $\Lambda$ passing close to
some point in $W^{s,+}_X(\sigma)$
(resp. $W^{s,-}_X(\sigma)$)
intersects $\Sigma^+$ (resp. $\Sigma^-$).
If $p\in \Lambda\cap\Sigma^+$ is close to $l^+$,
then
$W^s_X(p,\Sigma^+)$ is a vertical curve
crossing $\Sigma^+$.
If $p\in Per_X(\Lambda)$ and
$W^u_X(p)\cap W^{s,+}_X(\sigma)\neq\emptyset$, then
$W^u_X(p)$ contains an interval
$J=J_p$ intersecting $l^+$ transversally.
Similarly replacing $+$ by $-$.
\end{lemma}

To prove transitivity we shall use the
following criterium of Birkhoff.

\begin{lemma}
\label{trans}
Let $T$ be a compact, invariant set of $X$
such that for all open sets
$U,V$ intersecting $T$ there is
$s>0$ such that
$X_s(U\cap T)\cap V\neq\emptyset$.
Then $T$ is transitive.
\end{lemma}

Birkhoff's criterium is used together with
the following lemma.

\begin{lemma}
\label{*}
Let $\Lambda$ be a connected, singular-hyperbolic,
attracting set with dense periodic orbits
and only one singularity $\sigma$.
Let $U,V$ be open sets,
$p\in U\cap Per_X(\Lambda)$ and $q\in
V\cap Per_X(\Lambda)$.
If $W^u_X(p)\cap W^{s,+}_X(\sigma)\neq\emptyset$
and $W^u_X(q)\cap W^{s,+}_X(\sigma)\neq\emptyset$,
then there are $z\in W^u_X(p)$ arbitrarily close to $W^u_X(p)\cap W^{s,+}_X(\sigma)$
and $t>0$ such that
$X_t(z)\in V$.
A similar result holds replacing $+$ by $-$.
\end{lemma}

Hereafter we let $\Lambda$ be a
singular-hyperbolic set of
$X\in {\cal X}^r$
satisfying:

\begin{enumerate}
\item
$\Lambda$ is connected.
\item
$\Lambda$ is attracting.
\item
The closed orbits contained in $\Lambda$ are dense in
$\Lambda$.
\item
$\Lambda$ has only one singularity $\sigma$.
\end{enumerate}

We note that (3) above implies

\begin{description}
\item{(H1)}
$\Lambda=Cl(Per_X(\Lambda))$.
\end{description}

\begin{prop}
\label{co0}
Suppose that
the following property hold :
If $p,q\in Per_X(\Lambda)$ then
either
\begin{enumerate}
\item
$W^u_X(p)\cap W^{s,+}_X(\sigma)\neq\emptyset$
and $W^u_X(q)\cap W^{s,+}_X(\sigma)\neq\emptyset$ or
\item $W^u_X(p)\cap W^{s,-}_X(\sigma)\neq\emptyset$
and $W^u_X(q)\cap W^{s,-}_X(\sigma)\neq\emptyset$.
\end{enumerate}

Then, $\Lambda$ is transitive.
\end{prop}

\begin{proof}
By the Birkhoff's criterium
we only need to prove that
$\forall U,V$ open sets intersecting $\Lambda$
$\exists s>0$ such that
$X_s(U\cap \Lambda)\cap V\neq\emptyset$.
For this we proceed as follows:
By (H1) there are
$p\in Per_X(\Lambda)\cap U$
and $q\in Per_X(\Lambda)\cap V$.
First suppose that
the alternative (1)
holds. Then, by Lemma \ref{*},
there are $z\in W^u_X(p)$ and $t>0$ such that
$X_t(z)\in V$.
As $z\in W^u_X(p)$ we have that
$w=X_{-t'}(z)\in U$ for some $t'>0$.
As $\Lambda$ is an attracting set
we have that $w\in \Lambda$.
If $s=t+t'>0$ we conclude
that
$w\in (U\cap\Lambda)\cap X_{-s}(V)$ and so
$X_s(U\cap \Lambda)\cap V\neq\emptyset$.
If the alternative (2) of the corollary holds
we can find
$s>0$ such that $X_s(U\cap \Lambda)\cap V\neq\emptyset$
in a similar way (replacing $+$ by $-$).
The proof follows.
\end{proof}

\begin{prop}
\label{l2}
If there is a sequence $p_n\in Per_X(\Lambda)$ converging
to some point in $W^{s,+}_X(\sigma)$ such that
$
W^u_X(p_n)\cap W^{s,-}_X(\sigma)\neq\emptyset$
for all $n$, then
$\Lambda$ is transitive.
Similarly
interchanging the roles of $+$ and $-$.
\end{prop}

\begin{proof}
Let $p,q\in Per_X(\Lambda)$ be fixed.
Suppose that
$W^u_X(p)\cap W^{s,+}_X(\sigma)\neq\emptyset$ and $W^u_X(q)
\cap W^{u,-}_X(\sigma)\neq\emptyset$.
By Lemma \ref{sigma} we can fix a cross section $\Sigma=\Sigma^+$
through $W^{s,+}_X(\sigma)$
and a open arc $J\subset \Sigma\cap
W^u_X(p)$ intersecting $W^{s,+}_X(\sigma)$
transversally.
Again by Lemma \ref{sigma} we can assume that
$p_n\in \Sigma$
for every $n$.
Because
the direction $E^s$ of $\Lambda$
is contracting we have that
the size of $W^s_X(p_n)$ is uniformly
bounded away from zero.
It follows that there is $n$ large
so that $J$
intersects $W^s_X(p_n)$ transversally.
Applying the Inclination Lemma \cite{dMP} to the saturated
of $J\subset W^u_X(p)$,
and the assumption
$
W^u_X(p_n)\cap W^{s,-}_X(\sigma)\neq\emptyset
$,
we conclude that
$W^u_X(p)\cap W^{s,-}_X(\sigma)\neq\emptyset$.
So, the alternative (2) of Proposition
\ref{co0} holds.
It follows from this proposition that
$\Lambda$ is transitive.
\end{proof}

\begin{prop}
\label{co3}
If there is $a\in W^u_X(\sigma)\setminus \{\sigma\}$
such that
$\sigma\in \omega_X(a)$, then
$\Lambda$ is transitive.
\end{prop}

\begin{proof}
Without loss of generality we
can assume that there is
$z\in \omega_X(a)\cap W^+_X(\sigma)$.
If
$W^u_X(q)\cap W^{s,-}_X(\sigma)=\emptyset$
for all $q\in Per_X(\Lambda)$, then
$W^u_X(q)\cap W^{s,+}_X(\sigma)\neq\emptyset$
for all $q\in Per_X(\Lambda)$ (\cite{MPP}).
Then $\Lambda$ would be transitive by Proposition
\ref{co0} since the alternative
(1) holds $\forall p,q\in Per_X(\Lambda)$.
So, we can assume that there is
$q\in Per_X(\Lambda)$ such that
$W^u_X(q)\cap W^{s,-}_X(\sigma)\neq \emptyset$.
It follows from the dominating
condition of the singular-hyperbolic splitting of
$\Lambda$ that
the intersection $W^u_X(q)\cap W^{s,-}_X(\sigma)$
is transversal.
This allows us to choose
a point in
$W^u_X(q)$ arbitrarily close to
$W^{s,-}_X(\sigma)$ in the side accumulating
$a$.
As $\Lambda$ is attracting
and satisfies (H1),
it is not hard to find
a sequence $p_n\in Per_X(\Lambda)$
converging to
$z\in W^{s,+}_X(\sigma)$ such that for all $n$ there
is $p_n'$ in the orbit of $p_n$ such that the
sequence $p_n'$
converges to some point in $W^{s,-}_X(\sigma)$.
Now, suppose by contradiction that $\Lambda$ is not transitive.
Then Proposition \ref{l2}
would imply
$$
W^u_X(p_n')\cap W^{s,+}_X(\sigma)=\emptyset
\,\,\,\,\,\,\mbox{and}
\,\,\,\,\,\,\,W^u_X(p_n)\cap W^{s,-}_X(\sigma)=\emptyset
$$
for $n$ large.
But
$W^u_X(p_n)=W^u_X(p_n')$
since $p_n'$ and $p_n$ are in the same orbit of $X$.
So
$W^u_X(p_n)\cap (W^{s,+}_X(\sigma)\cup W^{s,-}_X(\sigma)
)=\emptyset$.
However
$$
W^u_X(p_n)\cap W^s_X(\sigma)=\emptyset,
$$
a contradiction
since $Sing_X(\Lambda)=\{\sigma\}$.
We conclude that $\Lambda$ is transitive
and the proof follows.
\end{proof}

\begin{thm}
\label{t1}
If $\Lambda$ is not transitive, then for all
$a\in W^u_X(\sigma)\setminus \{\sigma\}$
there is a periodic orbit
$O$ with positive expanding eigenvalues of $X$
such that $a\in W^s_X(O)$.
\end{thm}

\begin{proof}
We fix $a\in W^u_X(\sigma)\setminus \{\sigma\}$.
By contradiction we assume that
$\omega_X(a)$ is {\em not} a periodic orbit.
We obtain a contradiction once we prove that
if $p,q\in Per_X(\Lambda)$ then
$p,q$ satisfy one of the two alternatives
in Proposition \ref{co0}.
To prove this we proceed as follows.
As note before we have
that both
$W^u_X(p)$ and $W^u_X(q)$
intersect $W^s_X(\sigma)$ (\cite{MP1}).
Then we can assume

\begin{equation}
\label{estrela}
W^u_X(p)\cap W^{s,+}_X(\sigma)\neq\emptyset\,\,\,\,
\mbox{and}
\,\,\,\,
W^u_X(q)\cap W^{s,-}_X(\sigma)\neq\emptyset.
\end{equation}

By using (\ref{estrela})
and the linear coordinate around $\sigma$
it is easy to construct an open interval
$I=I_a$, contained in a suitable
cross-section $\Sigma=\Sigma_a$ of $X$ containing $a$,
such that $I\setminus \{a\}$ is formed by
two intervals $I^+\subset W^u_X(p)$ and
$I^-\subset W^u_X(q)$.
Observe
that the tangent vector of $I$ is contained
in $E^c\cap T\Sigma_a$.
Proposition \ref{co3} implies that
$\sigma\notin \omega_X(a)$
since
$\Lambda$ is not transitive.
It follows that
$H=\omega_X(a)$ is a hyperbolic
saddle-type set (\cite{MPP}).
As in \cite{M1} one proves
that $H$ is one-dimensional, and so,
the Bowen's Theory
of hyperbolic one-dimensional sets \cite{Bo}
can be applied.
In particular we can choose
a family
of cross-sections
${\cal S}=\{
S_1,\cdots ,S_r\}$ of small diameter
such that
$H$ is the flow-saturated
of $H\cap int({\cal S}')$, where
$
{\cal S}'=\cup S_i$
and $int({\cal S}')$ denotes the interior
of ${\cal S}'$.
On the other hand,
$I\subset
\Lambda$ since $\Lambda$
is attracting.
Recall that the tangent direction of $I$
is contained in $E^c$. Since
$E^c$ is volume expanding and
$H$ is non-singular we have that
the Poincar\'e map
induced by $X$ on ${\cal S}'$
is expanding along $I$.
As
in \cite[p. 371]{MP1}
we can find $\delta>0$ and a open arc sequence
$J_n\subset {\cal S}'$ in the positive orbit
of $I$ with lenght bounded away from $0$
such that there is
$a_n$ {\em in the positive orbit of $a$}
contained in the interior of
$J_n$.
We can fix $S=S_i\in {\cal S}$
in order to assume that
$J_n\subset S$ for
every $n$.
Let $x\in S$ be a limit point of
$a_n$.
Then $x\in H\cap int({\cal S}')$.
Because $I$ is tangent to
$E^c$ the interval sequence $J_n$
converges to an interval
$J\subset W^u_X(x)$ in the $C^1$ topology
($W^u_X(x)$ exists since $x\in H$ and $H$
is hyperbolic).
$J$ is not trivial since the lenght
of $J_n$ is bounded away from $0$.
If $a_n\in W^{s}_X(x)$
for $n$ large we would obtain that $x$
is periodic \cite[Lemma 5.6]{MP1},
a contradiction since $\omega_X(a)$ is not periodic.
We conclude that $a_n\notin W^s_X(x)$, $\forall n$.
As $J_n\to J$
and $\Lambda$ has strong stable manifolds
of uniformly size, we have that
$$
\exists z_n\in
(W^s_X(a_{n+1})\cap S)\cap
(J_{n}\setminus \{a_n\})
$$
for all $n$ large.
For every $n$ we
let $J_n^+$ and $J^-_n$ denote the two connected components
of $J_n\setminus \{a_n\}$ in a way that
$J^+_n$ is in the positive orbit of $I^+$ and
$J^-_n$ is in the positive orbit
of $I^-$. Clearly we have either $z_n\in J^+_n$ or $z_n\in J^-_n$.

If $z_n\in J^+_n$
there is
$v_{n+1}\in Per_X(\Lambda)\cap S$ close to $a_{n+1}$
such that
$$
W^s_X(v_{n+1})\cap J^+_{n}\neq\emptyset
\,\,\,\,
\mbox{and}
\,\,\,\,
W^s_X(v_{n+1})\cap J^-_{n+1}\neq\emptyset.
$$
Since $v_{n+1}$ is periodic
\cite{MPP} implies
that $W^u_X(v_{n+1})$ intersects
either $W^{s,+}_X(\sigma)$ or
$W^{s,-}_X(\sigma)$.
The choice of $v_{n+1}$ implies
that its orbits passes close to a point
in $W^{s,-}_X(\sigma)$.
Since $\Lambda$ is not transitive we conclude that
$W^u_X(v_{n+1})$ intersects
$W^{s,-}_X(\sigma)$
Since
$W^s_X(v_{n+1})\cap J^+_{n}$ is
transversal
the Inclination Lemma
then
imply $
W^u_X(p)\cap W^{s,-}_X(\sigma)\neq\emptyset$.
Hence
$$
W^u_X(p)\cap W^{s,-}_X(\sigma)\neq\emptyset\,\,\,\,
\mbox{and}
\,\,\,\,
W^u_X(q)\cap W^{s,-}_X(\sigma)\neq\emptyset.
$$

If $z_n\in J^-_n$ we can prove by similar arguments
that
$$
W^u_X(p)\cap W^{s,+}_X(\sigma)\neq\emptyset\,\,\,\,
\mbox{and}
\,\,\,\,
W^u_X(q)\cap W^{s,+}_X(\sigma)\neq\emptyset.
$$
These alternatives yield
the desired contradiction.
We conclude that
$\omega_X(a)=O$ for some
periodic orbit $O$ of $X$.
To finish let us prove that the expanding eigenvalue
of $O$ is positive.
Suppose by contradiction that it is not so.
Fix
a cross-section $\Sigma$
intersecting $O$ in a single point $p_0$.
This section defines a Poincar\'e map
$\Pi:Dom(\Pi)\subset \Sigma\to \Sigma$
of which $p_0$ is a hyperbolic fixed point.
The assumption implies that
$D\Pi(p_0)$ has negative expanding eigenvalue.
Because $p_0\in Per_X(\Lambda)$,
we have
that $W^u_X(p_0)$ intersects either
$W^{s,+}_X(\sigma)$ or
$W^{s,-}_X(\sigma)$ (\cite{MPP}).
We shall assume the first case
since the proof for the second one is similar.
We claim that
$W^u_X(p)\cap W^{s,+}_X(\sigma)\neq\emptyset$ for all
$p\in Per_X(\Lambda)$. Indeed,
let $p\in Per_X(\Lambda)$ be fixed.
Again
$W^u_X(p)$ intersects
$W^{s,+}_X(\sigma)$ or
$W^{s,-}_X(\sigma)$.
In the first case we are done. So, we can assume that
$W^u_X(p)\cap W^{s,-}_X(\sigma)\neq \emptyset$.
By flow-saturating this intersection
we obtain an interval $K\subset W^u_X(p)\cap \Sigma$
intersecting $W^s_X(p_0,\Sigma)$ transversally.
On the other hand, there is an interval
$J\subset W^{s,+}_X(\sigma)\cap \Sigma$
intersecting $W^u_X(p_0,\Sigma)$ transversally.
Since the expanding eigenvalue of $D\Pi(p_0)$
is negative we have by the Inclination Lemma
that
the backward iterates $\Pi^{-n}(J)$
of $J$ accumulates $W^s_X(p_0,\Sigma)$
{\em in both sides}.
Because $K$ has transversal intersection
with $W^s_X(p_0\Sigma)$ we conclude that
one of such backward iterates intersects
$K$.
This intersection point yields
$W^u_X(p)\cap W^{s,+}_X(\sigma)\neq \emptyset$ as desired.
This proves the claim.
The claim together with Proposition \ref{co0}
implies that $\Lambda$ is transitive, a contradiction.
This contradiction proves the result. 
\end{proof}

\begin{figure}[htv] 
\centerline{
\psfig{figure=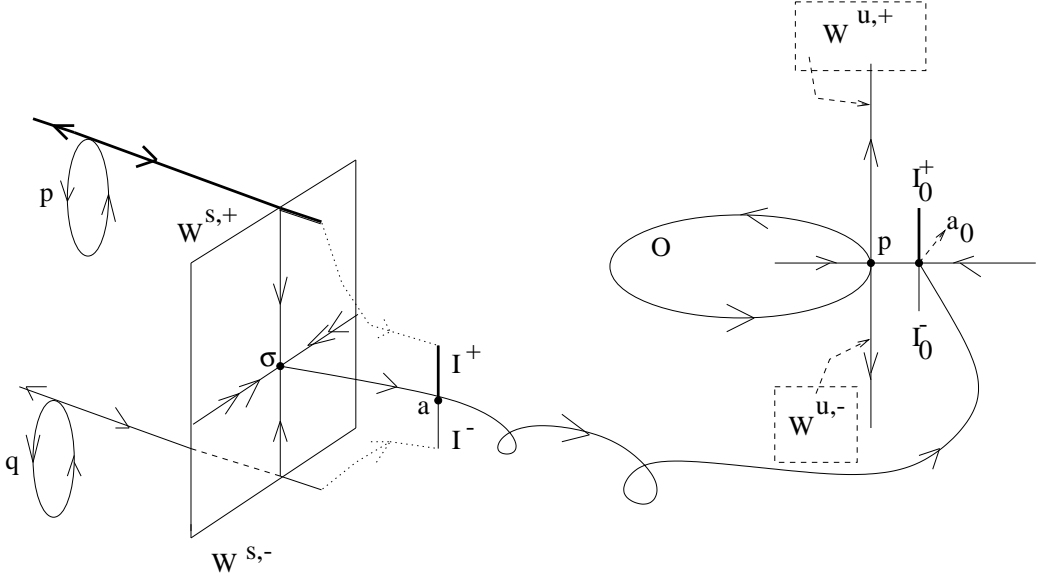,height=3in}}
\caption{\label{f.5} Definition of $W^{u,+}$ and
$W^{u,-}$.}
\end{figure}

Hereafter we shall assume that
$\Lambda$ is not transitive.
Let $a\in W^s_X(\sigma)\setminus \{\sigma\}$ be fixed.
By Theorem \ref{t1}
we have
that $a\in W^s_X(O)$ for some periodic orbit
$O$ with positive expanding eigenvalue.
This last property implies that
the unstable manifold
$W^u_X(O)$ of $O$ is a cylinder with generating curve
$O$. Then $O$ separates $W^u_X(O)$
in two connected components.
Such components are denoted by $W^{u,+},W^{u,-}$
according the following (see Figure \ref{f.5}):
There is
an interval
$I=I_a$, contained in a suitable cross-section of $X$ and containing $a$, such that
if $I^+,I^-$ are the connected components
of $I\setminus \{a\}$ then
$I^+\subset W^u_X(p)$ and $I^-\subset W^u_X(q)$ for some
periodic points $p,q\in \Lambda$
(recall that $\Lambda$ is not transitive).
In addition $I$ is tangent to the central direction $E^c$ of $\Lambda$ (see Figure \ref{f.5}).
As $a\in W^s_X(O)$ and $I$ is tangent to
$E^c$ we have that the flow
of $X$ carries $I$ to an interval $I'$
transverse to $W^s_X(O)$ at $a$.
Note that the flow carries $I^+$ and
$I^-$ into $I^+_0$ and $I^-_0$ respectively.

\begin{defi}
\label{u,+-}
We denote by
$W^{u,+}$ the connected component
of $W^{u}\setminus O$
which is accumulated
(via Inclination and Strong-$\lambda$-lemmas
\cite{dMP,D}) by the positive orbit
of $I^+_0$.
Similarly we denote
$W^{u,-}$ the connected component
of $W^{u}_X(O)\setminus O$
which is accumulated by the positive orbit
of $I^-_0$.
\end{defi}

The above definition does not depend on $p,q,J_p,J_q$
(this can be easily proved using
the Strong $\lambda$-lemma
\cite{D}).

\begin{prop}
\label{preliminar}
$W^{u,+}\cap W^{s,-}_X(\sigma)=\emptyset$
and $W^{u,+}\cap W^{s,+}_X(\sigma)\neq\emptyset$.
Similarly interchanging the roles of $+$ and $-$.
\end{prop}

\begin{proof}
For simplicity denote
$W=W^{u,+}$.
First we prove
$W\cap W^{s,-}_X(\sigma)=\emptyset$.
Suppose by contradiction that
$W\cap W^{s,-}_X(\sigma)\neq\emptyset$.
As this last intersection is transversal
there is an interval $J\subset
W^{s,-}_X(\sigma)$ intersecting
$W$ transversally.
Now, fix a cross-section $\Sigma=\Sigma^+$ as in
Lemma \ref{sigma} and let
$p\in Per_X(\Lambda)$ be such that
$W^u_X(p)\cap W^{s,+}_X(\sigma)\neq\emptyset$.
Then, there
is an small interval $I\subset
W^u_X(p)\cap \Sigma$ transversal to
$\Sigma\cap W^{s,+}_X(\sigma)$.
By the definition
of $W=W^{u,+}$
(Definition \ref{u,+-})
we have that the positive orbit
of $I$ accumulates
on $W$.
As $J$ is transversal to
$W$
the Inclination Lemma implies that
the positive orbit of $I$ intersects
$J$. This proves
$W^u_X(p)\cap W^{s,-}_X(\sigma)\neq\emptyset$
for all $p\in Per_X(\Lambda)$.
It would follow that alternative
(2) of Proposition \ref{co0}
holds $\forall p,q$ contradicting
the non-transitivity of $\Lambda$.
This contradiction proves
$W\cap W^{s,-}_X(\sigma)=\emptyset$
as desired.
Now suppose by contradiction that
$W\cap W^{s,+}_X(\sigma)=\emptyset$.
As $W\cap W^{s,-}_X(\sigma)=\emptyset$
we would obtain $W\cap W^s_X(\sigma)=\emptyset$
(\cite{MPP}).
But the denseness of the periodic orbits
together with the Inclination Lemma
imply $W\cap W^s_X(\sigma)\neq\emptyset$.
This is a contradiction which proves
$W\cap W^{s,+}_X(\sigma)\neq\emptyset$.
The result follows.
\end{proof}

\begin{prop}
\label{homoclinic3}
$H^+=Cl(W^{u,+})$ and $H^-=Cl(W^{u,-})$.
\end{prop}

\begin{proof}
Fix $q\in P^+$,
i.e.
$W^u_X(q)\cap W^{s,+}_X(\sigma)\neq\emptyset$.
Note that $W^{u,+}\cap
W^{s,+}_X(\sigma)\neq\emptyset$ by Lemma
\ref{preliminar}.
By using (H1) and the Inclination Lemma
it is not hard to prove that
$W^{u,+}$ accumulates on
$q$.
This proves
$
H^+\subset
Cl(W^{u,+})
$.
Conversely
let $x\in W^{u,+}$ be fixed.
By (H1) and
$W^{u,+}\subset \Lambda$
there is
$z\in Per_X(\Lambda)$ nearby $x$.
Choosing $z$ close to $x$
we assure
$W^s_X(z)\cap W^{u,+}\neq\emptyset$
because stable manifolds have size
uniformly bounded away from zero.
If $W^u_X(z)\cap W^{s,-}_X(\sigma)\neq \emptyset$ then
the Inclination Lemma
and
$W^s_X(z)\cap W^{u,+}\neq\emptyset$
would imply
$W^u\cap W^{s,-}_X(\sigma)\neq
\emptyset$.
This contradicts Proposition \ref{preliminar}
and so
$W^u_X(z)\cap W^{s,-}_X(\sigma)= \emptyset$.
By \cite{MP1} we obtain
$z\in P^+$ proving
$x\in H^+$.
The lemma is proved.
\end{proof}

\begin{prop}
\label{dense II}
If $z\in Per_X(\Lambda)$ and
$W^s_X(z)\cap W^{u,+}\neq\emptyset$,
then
$Cl(W^s_X(z)\cap W^{u,+})=
Cl(W^{u,+})$.
Similarly replacing $+$ by $-$.
\end{prop}

\begin{proof}
Previously we show
$Cl(W^{s,+}_X(\sigma)\cap W^{u,+})=Cl(W^{u,+})$.
Fix $x\in W^{u,+}$.
By (H1) there is
$w\in Per_X(\Lambda)$ close to $x$.
In particular $W^s_X(w)\cap W^{u,+}\neq\emptyset$.
If $W^u_X(w)\cap W^{s,+}_X(\sigma)
=\emptyset$
then $W^u_X(w)\cap W^{s,-}_X(\sigma)
\neq\emptyset$ by \cite{MP1}.
It would follow from
the Inclination Lemma that
$W^{u,+}\cap W^{s,-}_X(\sigma)\neq\emptyset$
contradicting Proposition \ref{preliminar}.
We conclude that 
$W^u_X(w)\cap W^{s,+}_X(\sigma)
\neq\emptyset$.
Note that $W^s_X(w)\cap W^{u,+}\neq \emptyset$
is close to $x$ as $w\to x$.
As
$W^u_X(w)\cap W^{s,+}_X(\sigma)
\neq\emptyset$ is transversal we
can apply the Inclination Lemma
in order to find
a transverse intersection
$W^{u,+}\cap W^{s,+}_X(\sigma)$ close to $x$. This proves
$Cl(W^{s,+}_X(\sigma)\cap W^{u,+})=Cl(W^{u,+})$.
Finally we prove
$Cl(W^s_X(z)\cap W^{u,+})=
Cl(W^{u,+})$.
Choose
$x\in W^{u,+}$.
As $Cl(W^{s,+}_X(\sigma)\cap W^{u,+})=Cl(W^{u,+})$
there is an interval
$I_x\subset W^{u,+}$ arbitrarily close
to $x$ such that
$I_x\cap W^{s,+}_X(\sigma)\neq\emptyset$.
The positive orbit
of $I_x$ first passes through
$a$ and, afterward, it accumulates
on $W^{u,+}$.
But $W^s_X(z)\cap W^{u,+}\neq\emptyset$ by assumption. As this intersection is transversal
the Inclination Lemma implies
that the positive orbit of
$I_x$ intersects $W^s_X(z)$.
By taking the backward flow of the last intersection
we get $W^s_X(z)\cap I_x\neq\emptyset$.
This proves the lemma.
\end{proof}

Given $z\in Per_X(\Lambda)$ we denote by
$H_X(z)$ the homoclinic class associated to $z$.

\begin{prop}
\label{cl=homoclinic}
If $z\in Per_X(\Lambda)$ is close to a point in
$W^{u,+}$, then
$H_X(z)=Cl(W^{u,+})$.
Similarly replacing $+$ by $-$.
\end{prop}

\begin{proof}
Let $z\in Per_X(\Lambda)$ be a point
close to one in $W^{u,+}$.
It follows from the continuity
of the stable manifolds that
$
W^s_X(z)\cap W^{u,+}\neq\emptyset$.
We claim that
$H_X(z)=Cl(W^{u,+})$.
Indeed, by Proposition \ref{preliminar}
one has $W^{u,+}\cap W^{s,-}_X(\sigma)=
\emptyset$. This equality and the Inclination Lemma
imply $W^u_X(z)\cap W^s_X(\sigma)\subset
W^{s,+}_X(\sigma)$.
By
Proposition \ref{dense II} we have that
$W^s_X(z)\cap W^{u,+}$ is dense
in $W^{u,+}$ since
$W^s_X(z)\cap W^{u,+}\neq\emptyset$.
Let $\Sigma$ be a cross-section containing $p_0$
and
fix $x\in W^{u,+}$.
We can assume $x,z\in \Sigma$.
Since $W^u_X(z)\cap W^s_X(\sigma)\neq\emptyset$
and
$W^u_X(z)\cap W^s_X(\sigma)\subset W^{s,+}_X(\sigma)$
there is an interval
$I\subset W^u_X(z)$ intersecting
$W^{s,+}_X(\sigma)$.
Then, the positive orbit of $I$ yields the interval
$J$ in the figure.
In addition,
the positive orbit of $J$ yields the interval
$K$ in that figure.
Note that the positive orbit of $K$
{\em accumulates} $W^{u,+}$ (recall Definition
\ref{u,+-}).
As $W^s_X(z)\cap W^{u,+}$
is dense in $W^{u,+}$
and $x\in W^{u,+}$
we have that
$W^s_X(z)$ passes close to $x$ as indicated in the figure.
The Inclination Lemma applied to the positive orbit
of $K$ yields a homoclinic point
$z'$ associated to $z$ which is close to $x$.
This proves that $x\in H_X(z)$, and so, $Cl(W^{u,+})\subset H_X(z)$.
The reversed inclusion
is a direct consequence of the Inclination Lemma
applied to $W^s_X(z)\cap W^{u,+}\neq\emptyset$.
We conclude that $Cl(W^{u,+})=
H_X(z)$ proving the result.
\end{proof}

\begin{thm}
\label{homoclinic class}
Let $\Lambda$ be a singular-hyperbolic set
of a $C^r$ flow $X$ on a closed three-manifold,
$r\geq 1$.
Suppose that the properties below hold.

\begin{enumerate}
\item
$\Lambda$ is connected.
\item
$\Lambda$ is attracting.
\item
The closed orbits contained in $\Lambda$ are dense in $\Lambda$.
\item
$\Lambda$ has a unique singularity $\sigma$.
\item
$\Lambda$ is not transitive.
\end{enumerate}
Then,
$H^+$
and $H^-$
are homoclinic classes of $X$.
\end{thm}

\begin{proof}
Let $\Lambda$ be a singular-hyperbolic set
of $X$ satisfying (1)-(5) of
Theorem \ref{homoclinic class}.
Then we can apply the results in this section.
To prove that $H^+$ is
a homoclinic classs it suffices
by Proposition \ref{homoclinic3} to
prove that $Cl(W^{u,+})$ is a homoclinic
class.
By (3) of Theorem \ref{homoclinic class}
we can choose $z\in Per_X(\Lambda)$
arbitrarily close to
a point in $W^{u,+}$.
Then
$Cl(W^{u,+})=H_X(z)$
by Proposition \ref{cl=homoclinic}
and the result follows.
\end{proof}

\section{Proof of Theorem \ref{thA}}

First we introduce some notations.
Hereafter $M$ is a compact $3$-manifold
and ${\cal X}^r$ is
the space of $C^r$ flows in $M$ equiped
with the $C^r$ topology, $r\geq 1$.
The {\em nonwandering set}
of $X\in{\cal X}^r$ is the set
$\Omega(X)$ of points $p\in M$ such that
for all neighborhood
$U$ of $p$ and $T>0$ there is $t>T$ such that
$X_t(U)\cap U\neq\emptyset$.
An attracting $\Lambda$ with isolating
block $U$ has a continuation
$\Lambda(Y)$ for $Y$ $C^r$ close to $X$ defined
by
$
\Lambda(Y)=\cap_{t>0}Y_t(U)
$. This continuation is
then definied when $\Lambda$ is an attractor.
A compact invariant set is {\em non-trivial}
if it is not a closed orbit of $X$.
Transitive sets for flows are always connected.
The proof of Theorem \ref{thA} is based on the
following result.

\begin{thm}
\label{thA'}
Let $\Lambda$ be a singular-hyperbolic set
of $X\in {\cal X}^r$,
$r\geq 1$.
Suppose that the following properties hold.

\begin{enumerate}
\item
$\Lambda$ is connected.
\item
$\Lambda$ is attracting.
\item
The closed orbits contained in $\Lambda$ are dense in
$\Lambda$.
\item
$\Lambda$ has only one singularity.
\item
$\Lambda$ is not transitive.
\end{enumerate}
Then, for every neighborhood $U$
of $\Lambda$ there is a flow $Y$
$C^r$ close to $X$ such that
$$
\Lambda(Y)\not \subset\Omega(Y).
$$
\end{thm}

To prove this theorem we shall
use the following definitions and facts.
Let $X\in {\cal X}^r$ and
$\Lambda$ be a singular-hyperbolic
set of $X$
satisfying (1)-(5) of Theorem \ref{thA'}.
Let $\sigma$ be the unique singularity of $\Lambda$.
As mentioned in the previous section $\sigma$ is Lorenz-like.
As in Section 2 we have that $W^{ss}_X(\sigma)$ divides
$W^s_X(\sigma)$ in two connected components
which we denote by $W^{s,+}_X(\sigma),W^{s,-}_X(\sigma)$.
For simplicity we denote such components
by $W^{s,+},W^{s,-}$ respectively.
Recall that $Per_X(\Lambda)$ denotes
the union of the
periodic orbits of $X$ in $\Lambda$.
A point in $Per_X(\Lambda)$ is said to be periodic.
Fix such $a\in W^u_X(\sigma)\setminus \{\sigma\}$.
By Theorem \ref{t1} we have that
$\omega_X(a)=O$ for some periodic orbit
with positive expanding eigenvalues of $X$.
In particular, $W^{u,+}$ and $W^{u,-}$ are defined
(Definition \ref{u,+-}).

\begin{lemma}
\label{l}
$Cl(W^{u,+})\cap W^{s,-}=\emptyset$.
\end{lemma}

\begin{proof}
Suppose by contradiction that
$Cl(W^{u,+})\cap W^{s,-}\neq\emptyset$.
Choose a cross-section
$\Sigma^-$ intersecting $W^{s,-}_X(\sigma)$
such that every orbit's sequence in
$\Lambda$ converging to some point
in $W^{s,-}$ intersect $\Sigma^-$
(Lemma \ref{sigma}).
Let $q\in \Lambda$ be periodic
such that $W^u_X(q)\cap W^{s,-}_X(\sigma)
\neq\emptyset$.
As $Cl(W^{u,+})\cap W^{s,-}\neq\emptyset$
we have that $Cl(W^{u,+})\cap \Sigma^-\neq\emptyset$.
Because the closed orbits
are dense
we can prove that $q\in Cl(W^{u,+})$.
It follows that $H^-
\subset Cl(W^{u,+})$, and so,
$\Lambda=Cl(W^{u,+})$ by Lemma \ref{cara}.
On the other hand, since $\Lambda$ is not transitive,
we have that
$H^+=Cl(W^{u,+})$ (Proposition \ref{homoclinic3})
and $H^+$ is a homoclinic class
(Theorem \ref{homoclinic class}).
Since homoclinic classes are transitive sets
we conclude that $\Lambda$ is transitive, a contradiction. This contradiction proves the result.
\end{proof}

\begin{lemma}
\label{l2'}
Let $D$ be a fundamental domain
of $W^{uu}_X(p_0)$ contained in $W^{u,+}$.
Then there are a neighborhood
$V$ of $D$ and a cross-section
$\Sigma^-$ of $X$
intersecting $W^{s,-}$ satisfying the following properties.

\begin{enumerate}
\item
Every $X$-orbit's sequence in
$\Lambda$ converging to some point
in $W^{s,-}$ intersects $\Sigma^-$.
\item
Every positive $X$-orbit with
initial point in
$V$ does not interset $\Sigma^-$.
\end{enumerate}
\end{lemma}

\begin{proof}

Fix a fundamental domain
$F$ of $W^s_X(\sigma)$ and define
$F^-=F\cap W^{s,-}_X(\sigma)$
Then there is a compact interval
$F'\subset F^-$
such that $\Lambda\cap F^-\subset F'$.
By Lemma \ref{l} there is
$\epsilon>0$ such that
$B_\epsilon(Cl(W^{u,+}))\cap B_\epsilon(F')=\emptyset$.
Clearly we can choose a cross-section
$\Sigma^-$ of $X$ inside $B_\epsilon(F')$
such that
every $X$-orbit's sequence in
$\Lambda$ converging to some point
in $W^{s,-}$ intersects $\Sigma^-$.
On the other hand,
as $W^{u,+}$ is invariant
and $Cl(W^{u,+})\cap B_\epsilon(F')=\emptyset$
we have that
every positive orbit
with initial point
in $D$ cannot intersect $\Sigma^-$.
By using the contracting foliation
of $\Lambda$ we have the same property
for every positive trajectory with
initial point in a neighborhood
$V$ of $D$.
This proves the result.
\end{proof}

Now we define a perturbation
(pushing) close to a point $a\in
W^u_X(\sigma)\setminus \{\sigma\}$.
For this end we fix the following cross-sections:

\begin{enumerate}
\item
$\Sigma_a$
containing $a$ in its interior.
\item
$\Sigma'=X_1(\Sigma)$.
\item
$\Sigma_0$
intersecting $O$ in a single interior point.
\item
$\Sigma^+,\Sigma^-$ which
intersect $W^{s,+},W^{s,-}$
and pointing to the side of $a$ respectively.
\end{enumerate}

\begin{figure}[htv] 
\centerline{
\psfig{figure=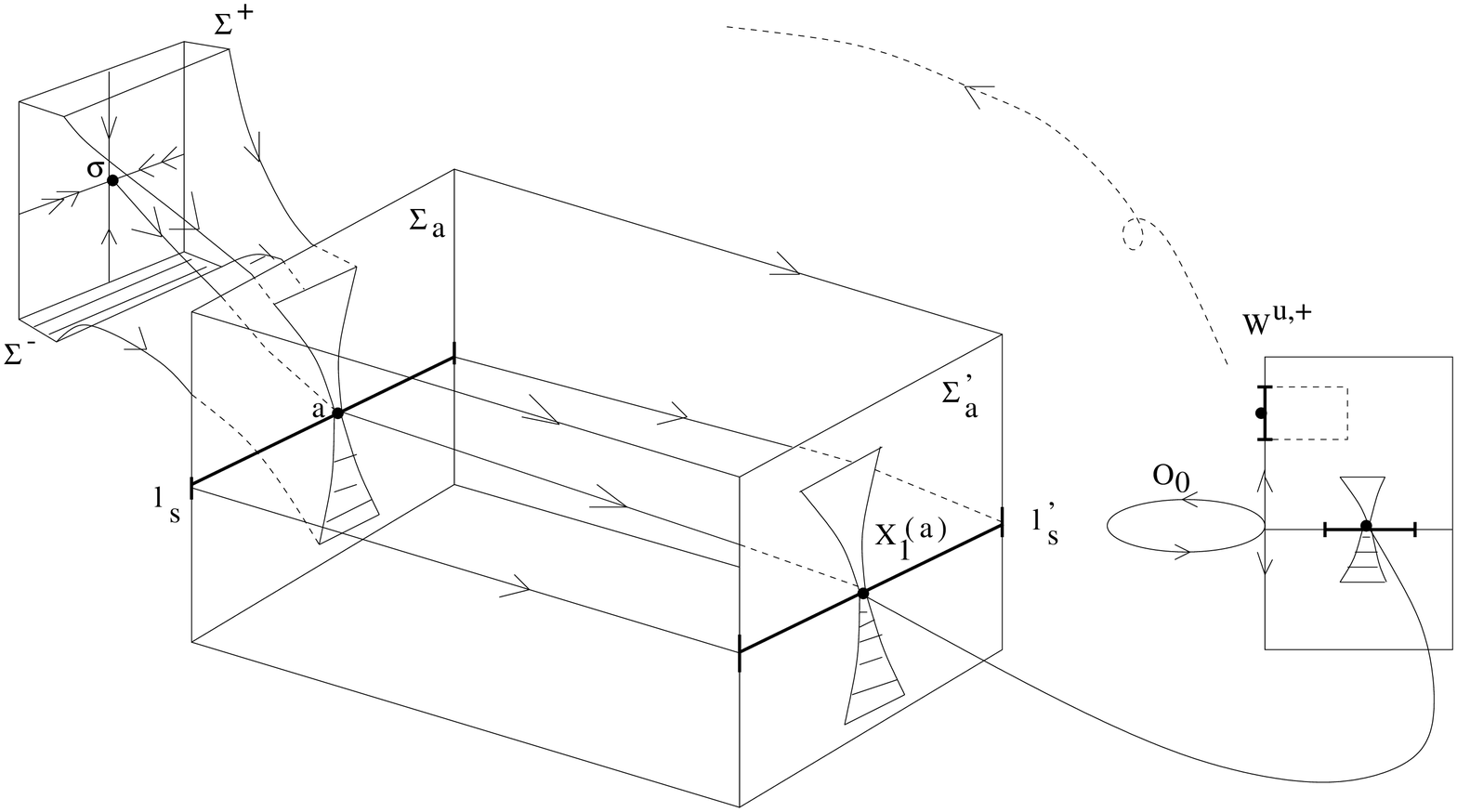,height=2.5in}}
\caption{\label{f.1} Unperturbed flow $X$.}
\end{figure}

\begin{figure}[htv] 
\centerline{
\psfig{figure=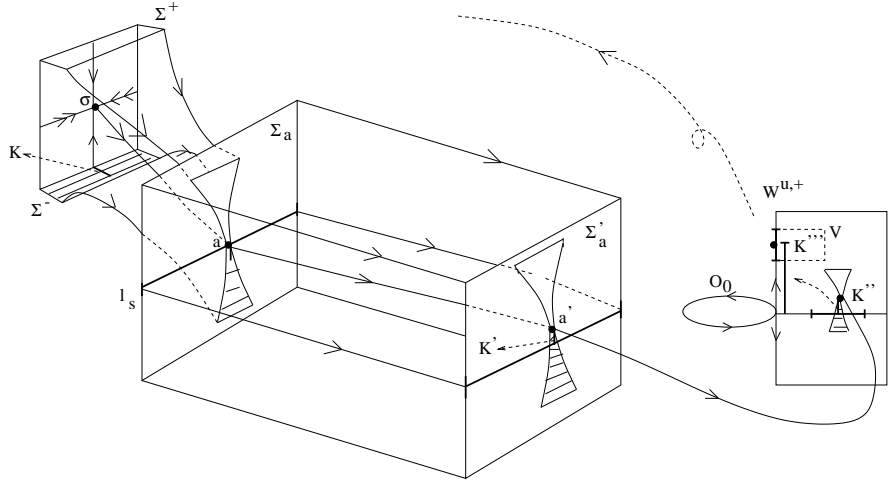,height=2.5in}}
\caption{\label{f.2'} Perturbed flow $Y$.}
\end{figure}

Every $X$-orbit intersecting
$\Sigma^+\cup \Sigma^-$ will
intersect $\Sigma$.
Note that there is a well defined neighborhood
${\cal O}$
given by
$$
{\cal O}=
X_{[0,1]}(\Sigma_a).
$$
This neighborhood will be the support
of the pushing described
in the Figures \ref{f.1} and \ref{f.2'}.
The pushing in
${\cal O}$ yielding the
perturbed flow $Y$ of $X$
is obtained in the standard
way (see \cite{dMP}).

We have to prove that
$\Lambda(Y)\not \subset \Omega(Y)$ for the perturbed
flow $Y$.
For this purpose we observe that
by (5) of Theorem \ref{thA'}
and Proposition \ref{co0}
we can assume that there
$q$ periodic in
$U$ such that
$W^u_X(q)\cap W^{s,-}_X(\sigma)\neq\emptyset$.
We obtain in this way an interval
$K$ in $\Sigma^-\cap W^u_X(q)$
crossing $\Sigma^-$ as in Figure \ref{f.2'}.
The $Y$-flow carries $K$ to an interval
$K''$ as in Figure \ref{f.2'}.
Let $q(Y),W^u_Y(q(Y)),K''(Y),\sigma(Y)$ denote
the continuation of these objects for
the perturbed flow $Y$.
We observe that $K(Y)\subset \Lambda(Y)$ since
$\Lambda(Y)$ is an attracting set,
$q(Y)\in \Lambda(Y)$ and $K\subset W^u_Y(q(Y))$.
Then Theorem \ref{thA'} will follow from
the lemma below.

\begin{lemma}
\label{lemma*}
$K(Y)\not \subset \Omega(Y)$.
\end{lemma}

\begin{proof}
Suppose by contradiction that
$K(Y)\subset \Omega(Y)$ and
pick $x\in Int(K(Y))$ ($=$ the interior of the interval
$K(Y)$).
The $Y$'s flow carries the points
nearby $x$ to the neighborhood
$V$ depicted in Figure \ref{f.2'}.
This neighborhood is
obtained by saturation
a fundamental domain in $W^{u,+}$
by the strong stable manifolds \cite{HPS}.
Note that there are points $x'$ nearby $x$ which
back close to $x$ by the positive $Y$-flow
($x\in K(Y)\subset \Omega(Y)$).
In particular,
{\em the positive $Y$-orbit of $x'$ returns to $\Sigma^-$}.
On the other hand, Lemma \ref{l2'}-(2) implies that
every $X$-orbit starting in $V$
doesn't intersect $\Sigma^-$.
Since $X=Y$ outside ${\cal O}$
we conclude that the positive
$Y$-orbit of $x'$ intersects $\Sigma^+$.
Afterward this positive orbit
passes throught the box ${\cal O}$
and arrives to $V$. By repeating the argument
we conclude that {\em the positive $Y$-orbit
of $x'$
never return to $\Sigma^-$}, a contradiction.
The lemma is proved.
\end{proof}

{\flushleft{\bf Proof of Theorem \ref{thA}: }}
Let $\Lambda$ be a singular-hyperbolic attractor
of a $C^r$ flow $X$ on a compact $3$-manifold $M$.
Assume that $\Lambda$
is $C^r$ critically-robust and
has a unique singularity $\sigma$.
Denote
by $\Lambda(Y)=\cap_{t>0}Y_t(U)$
the continuation of $\Lambda$
in a neighborhood $U$ of $\Lambda$
for $Y$ close to $X$.
Denote by $C(Y)$ the union of the closed orbits
of a flow $Y$.
As $\Lambda$ is $C^r$ critically-robust
it follows that
there is a neighborhood $U$
of $\Lambda$ such that
$\Lambda(Y)\cap C(Y)$ is dense in $\Lambda(Y)$
for every flow $Y$ $C^r$ close to $X$.
Clearly $\Lambda(Y)$ is a singular-hyperbolic set
of $Y$ for all $Y$ close to $X$.
Because $\Lambda$ has a unique singularity
we have that $\Lambda(Y)$ has a unique
singularity as well.
Because $\Lambda$ is transitive we have that $\Lambda$ is connected.
Then the neighborhood
$U$ above can be arranged connected, and so,
$\Lambda(Y)$ is also connected.
Summarizing we have that
$\Lambda(Y)$ is a singular-hyperbolic set of $Y$
satisfying (1)-(4) of Theorem \ref{thA'}.
If $\Lambda$ were not $C^r$ robust, then
it would exist $Y$ $C^r$ close to $X$ such that
$\Lambda(Y)$ is not transitive.
In this case $\Lambda(Y)$ satisfies
(1)-(5) of Theorem \ref{thA'}. It follows that
there is $Y'$ $C^r$ close to $Y$
such that
$\Lambda(Y')\not \subset \Omega(Y')$. And this is
a contradiction since
$\Lambda(Y')
\subset \Omega(Y')$
(as $\Lambda(Y')\cap C(Y')$ is dense in $\Lambda(Y')$).
This proves the result.
\qed

\vspace{0.2cm}
\noindent C. A. Morales, M. J. Pacifico\\
Instituto de Matem\'atica\\
Universidade Federal do Rio de Janeiro\\
C. P. 68.530, CEP 21.945-970\\
Rio de Janeiro, R. J. , Brazil\\
e-mail: {\em morales@impa.br}, {\em pacifico@impa.br}

\end{document}